\documentclass{amsart}
\usepackage{amsmath}
\usepackage{amssymb}
\usepackage[english]{babel}
\usepackage{graphics}
\usepackage[all]{xy}
\usepackage{amsrefs}
\usepackage{amsmath}
\usepackage{amsthm}
\usepackage{amssymb}
\usepackage{amsfonts}
\usepackage{amsxtra}
\usepackage{tikz}
\usepackage{multicol}
\usepackage{subfigure}
\usepackage{centernot}
\usepackage{soul}
\usepackage{xcolor}
\usepackage{layout}
\usepackage{geometry}
\usepackage{amsthm}

\usepackage{amstext} 
\usepackage{array}   
\newcolumntype{L}{>{$}c<{$}} 
\usepackage{longtable}

\usepackage[colorlinks=true,linkcolor=blue,urlcolor=blue, citecolor=blue]{hyperref}
\usepackage{fullpage}
\usepackage{epsfig}
\usepackage{verbatim}
\usepackage{enumitem}
\usepackage{algorithmicx}
\usepackage{float}
\usepackage{algpseudocode}

\usepackage[justification=centering]{caption}

\newcommand{\divides}{\mid}

\newtheorem{Theorem}{Theorem}[section]

\newtheorem{Lemma}[Theorem]{Lemma}

\theoremstyle{definition}

\theoremstyle{remark}
\newtheorem*{Remark}{Remark}

\usepackage{graphicx}
\usepackage{adjustbox}

\title{Classifying Solvable Primitive Permutation Groups of Low Rank}

\author[M.~Dolorfino]{Mallory~Dolorfino}
\author[L. ~Martin]{Luke~Martin}
\author[Z. ~Slonim]{Zachary Slonim}
\author[Y. ~Sun]{Yuxuan Sun}
\author[Y. ~Yang]{Yong Yang}
\address{Mallory~Dolorfino\\ Kalamazoo College \\ Kalamazoo, Michigan, USA \\
  \href{mailto:mallory.dolorfino19@kzoo.edu}
  {{\ttfamily\upshape mallory.dolorfino19@kzoo.edu}}}
\address{Luke~Martin\\ Gonzaga University\\ Spokane, Washington, USA \\
  \href{mailto:lwmartin2019@gmail.com}
  {{\ttfamily\upshape lwmartin2019@gmail.com}}}
\address{Zachary~Slonim\\ University of California, Berkeley \\ Berkeley, California, USA \\
  \href{mailto:zachslonim@berkeley.edu}
  {{\ttfamily\upshape zachslonim@berkeley.edu}}}
\address{Yuxuan~Sun\\ Haverford College \\ Haverford, Pennsylvania, USA \\
  \href{mailto:ysun1@haverford.edu}
  {{\ttfamily\upshape ysun1@haverford.edu}}}
\address{Yong~Yang\\ Texas State University \\ San Marcos, Texas, USA \\
  \href{mailto:yang@txstate.edu}
  {{\ttfamily\upshape yang@txstate.edu}}}

\begin{document}

\maketitle

\begin{abstract}

    Suppose that a finite solvable permutation group $G$ acts faithfully and primitively on a finite set $\Omega$. Let $G_0$ be the stabilizer of a point $\alpha \in \Omega$ and the rank of $G$ be the number of distinct orbits of $G_0$ in $\Omega$ (including the trivial orbit $\{\alpha\}$). Then $G$ always has rank greater than four except for in a few cases. We completely classify these cases in this paper.

\end{abstract}

\begin{section}{Introduction}

Suppose that $G$ is a finite, transitive, solvable permutation group acting on a set $S$ with $n$ elements. Let $G_0$ be the stabilizer of a point $\alpha \in \Omega$. Define the rank of a permutation group, denoted $r(G),$ as the number of distinct orbits of $G_0$ in $S$ (including the trivial orbit $\{\alpha\}$). Huppert \cite{Huppert} and Foulser \cite{Foulser} classified all finite, solvable, permutation groups of rank two and three respectively, and Foulser restricted the rank four groups to a small list of possibilities. This paper completes the classification of all groups of rank less than $5$ by explicitly confirming these past results and computationally constructing the groups of rank $4$.

Let $G$ be a primitive, solvable, permutation group acting on a set with $n$ elements. Then $n = p^k$ for some prime $p$ and positive integer $n$. Also, $G$ must include a unique, normal minimal subgroup $M$, where $M$ is an elementary abelian group such that $|M| = p^d$. Thus $M$ behaves as a $d$-dimensional vector space over $\mathrm{GF}(p),$ and we use the notation $V :=M$ \cite{Seager}*{Section 2}.

We also have that $G = V \rtimes G_0$, where $G_0$ acts on $V$ as an irreducible subgroup of $\mathrm{GL}(V)$. The converse is true as well. That is, if we have an irreducible subgroup, $\overline{G}$ of $\mathrm{GL}(V)$, we can construct a primitive permutation group by the semidirect product $G = V \rtimes \overline{G}$ \cite{Seager}*{Section 2}. Thus, all the rank 4 cases Foulser suggests can be constructed by considering a solvable linear group acting on a vector space $V$, giving a primitive, solvable, permutation group.

We follow Foulser's analysis of solvable, linear groups which divides this class of groups into three subclasses, $\mathfrak{A}$, $\mathfrak{B}$, and $\mathfrak{L}$. $\mathfrak{A}$ consists of subgroups of $\Gamma (V)$, the semilinear group of $V$, $\mathfrak{B}$ consists of the remaining primitive linear groups, and $\mathfrak{L}$ consists of those which are imprimitive. When $G_0 \in \mathfrak{B}$, we further divide $\mathfrak{B}$ into two subclasses.

Let $G_0 \in \mathfrak{B}$ be a finite, solvable group which acts faithfully, irreducibly, and primitively on a finite vector space $V$ of dimension $d$ over a finite field $\mathbb{F}$ of characteristic $p$. From Foulser's work, $G_0$ contains two normal subgroups, $U \leq F$ such that $|F : U| = q^{2m}$ for some prime $q.$ Then Theorem \ref{thm21} below guarantees that $G_0$ has a unique extraspecial subgroup $E$ of order $q^{2m+1}$ where $|E/Z(E)|=q^{2m}$ is a prime power.

$\mathfrak{A}$ contains the groups whose extraspecial subgroups $E$ are trivial, while  $\mathfrak{B}_I$ and $\mathfrak{B}_R$ contain those whose extraspecial subgroups $E$ that act irreducibly and reducibly on $V$, respectively.

In this paper, we explicitly construct these groups in M\textsc{agma} \cite{MAGMA} and G\textsc{ap} \cite{GAP}, confirming the results of Foulser for rank 3 groups and fully classifying the rank 4 groups. With this classification, we expect that many proofs of past results could be simplified, and these constructions will have future applications.

Computationally, we prove the following result.

\begin{Theorem}\label{thm11}
    Let $G$ be a finite, solvable, primitive, permutation group. If $r(G) \leq 4,$ then (at least) one of the following holds:
    \begin{enumerate}
        \item $G_0 \in \mathfrak{A};$
        \item $G_0 \in \mathfrak{B}_I$ and $G_0$ is listed in Table \ref{tbl41};
        \item $G_0 \in \mathfrak{B}_R$ and $G_0$ is listed in Table \ref{tbl42};
        \item $G_0 \in \mathfrak{L}$ and there exists an imprimitivity decomposition $V=\sum_{i=1}^{r} \oplus V_{i}$ for $r=2$ or 3, and $G_{0} \mid V_{i}$ is transitive on $V_{i}-\{0\}$, for $1 \leq i \leq r$.
    \end{enumerate}
\end{Theorem}

\end{section}

\begin{section}{Theoretical Background}\label{theory}

In this section, we present the theoretical background motivating our work. First note that if $G_0$ acts primitively then it acts quasi-primitively, which means that all non-trivial normal subgroups of $G_0$ act homogeneously on $V$. We now describe the structure of a finite, solvable group $G_0$ that acts faithfully, irreducibly and quasi-primitively on a finite vector space $V$ of dimension $d$ over a finite field $\mathbb{F}$ of characteristic $p$.

\begin{Theorem}\label{thm21}
    \cite{Yang}*{Theorem 2.1, 2.2} Suppose that a finite solvable group $G_0$ acts faithfully, irreducibly, and quasi-primitively on a $d$-dimensional finite vector space $V$ over a finite field $\mathbb{F}$ of characteristic $p.$ Then every normal abelian subgroup of $G_0$ is cyclic, and $G_0$ has normal subgroups $Z \leq U \leq F \leq A \leq G_0$ and a characteristic subgroup $E \leq F.$ Suppose $|F : U| := q^{2m}$ is a prime power. Then the following statements hold:

    \begin{enumerate}
        \item $F = EU$ is a central product where $Z := E \cap U = \mathbf{Z}(E)$ and $\mathbf{C}_{G_0}(F) \leq F;$
        \item $F/U \cong E/Z$ is a direct sum of completely reducible $G_0/F$-modules;
        \item $E$ is an extraspecial $q-group$ and $|E/Z| = q^{2m}.$ Furthermore, $q^m$ divides $d$ and $q \neq p;$
        \item $A = \mathbf{C}_{G_0}(U),$ $G_0/A \lesssim \mathrm{Aut}(U),$ and $A/F$ acts faithfully on $E/Z;$
        \item $A/\mathbf{C}_A(E/Z) \lesssim \mathrm{Sp}(2m,q);$
        \item $U$ is cyclic and acts fixed point freely on $W,$ where $W$ is an irreducible submodule of $V_{U};$
        \item $|U|$ divides $p^k - 1$ for some $k \geq 1,$ and $W$ can be identified with the span of $U$ which is isomorphic to $\mathrm{GF}(p^k);$
        \item $|V| := n = |W|^{q^mb}$ for some integer $b;$
        \item $G_0/A$ is cyclic and $|G_0:A| \divides \mathrm{dim}(W).$ Further, $G_0 = A$ when $q^m = d.$
    \end{enumerate}
\end{Theorem}

The next theorem, due to Fousler, gives the partial classification of all finite, solvable, permutation groups $G$ of low rank $(\leq 4)$ by partitioning $G_0$ into three classes.

\begin{Theorem}\label{foulserthm}
     \cite{Foulser}*{Theorem 1.2} Let $G$ be a finite solvable primitive permutation group of degree $n.$ Then $r(G) \geq 5,$ except possibly in the following cases:

    \begin{enumerate}
        \item $G_0 \in \mathfrak{A};$
        \item $G_0 \in \mathfrak{B}_I$ and one of the following cases applies:
            \begin{enumerate}
                \item $q^m = 3, p^k = 4$ or $7,$ and $n = p^{3k};$
                \item $q^m = 2, p^k \leq 71,$ and $n=p^{4k};$
                \item $q^m = 4, p^k = 3, 5,$ or $7,$ and $n = p^{4k}.$
            \end{enumerate}
        \item $G_0 \in \mathfrak{B}_R$ and one of the following cases applies:
            \begin{enumerate}
                \item $q^m = 2, p^k \leq 7,$ and $n = p^{4k};$
                \item $q^m = 2, p^k = 3,$ and $n = 3^6;$
                \item $q^m = 2, p^k = 3,$ and $n = 3^{10}.$
            \end{enumerate}
        \item $G_0 \in \mathfrak{L}$, there exists an imprimitivity decomposition $V=\sum_{i=1}^{r} \oplus V_{i}$ for $r=2$ or 3, and $G_{0} \mid V_{i}$ is transitive on $V_{i}-\{0\}$, for $1 \leq i \leq r$.
    \end{enumerate}
\end{Theorem}

\begin{Remark}
    All of the groups of rank $3$ are fully classified by Foulser \cite{Foulser}*{Theorem 1.1}.
\end{Remark}

We now present the structure of the groups listed in Table \ref{tbl41} and Table \ref{tbl42}. We follow the structure given in Theorem \ref{thm21} using a layered approach, starting from the extraspecial group and working our way up.

We start by outlining some properties of the extraspecial subgroup $E$. Extraspecial groups, $E$, are special types of $p$-groups where $\mathbf{Z}(E)$ is a cyclic group of order $p$ and $E / \mathbf{Z}(E)$ is a non-trivial, elementary abelian $p$-group. Namely, $E$ cannot be abelian. Any extraspecial group $E$ has order $q^{2m+1}$ for some prime $q$ and positive integer $m.$ Furthermore, there are exactly two isomorphism classes of extraspecial groups of order $q^{2m+1}$ for any given $q$ and $m.$ Both of these are isomorphic to central products of $m$ extraspecial groups of order $q^3$ \cite{Gorenstein}*{Chapter 5, Theorem 5.2(ii)}.

Given the restrictions Foulser places on the possible low rank groups, we only have three distinct cases to consider: $q^m = 2^1, 3^1,$ and $2^2$. For $q = 2$, the two extraspecial groups of order $q^3$ are $D_8$ and $Q_8$, the dihedral and quaternion groups of order $8$. Both of these groups have exponent exp$(E) = 4$. For $q$ odd, the two extraspecial groups of order $q^3$ are isomorphic to either the three-dimensional uni-triangular matrices with entries in $\mathbb{F}_q$ or the group $\mathbb{Z}_q \rtimes \mathbb{Z}_{q^2}$. These groups have exponent $q$ and $q^2,$ respectively. Finally, for $q=2, m=2$, the extraspecial groups of order $2^5$ are the central products $D_8 * D_8$ and $D_8 * Q_8$, which we will denote by $E^+$ and $E^-$ respectively. This last statement follows from the fact that $D_8 * D_8 \cong Q_8 * Q_8$. Again, both of these groups have exponent $4$. We now make a claim which will further restrict the possibilities for $E$ in the first case.

\begin{Lemma}\label{extraspeciallimits}

    Let $G$ be a finite, solvable, primitive, permutation group of degree $n$ whose subgroup $G_0$ belongs to either case $(2)$ or $(3)$ in Theorem \ref{foulserthm} and $E$ be its extraspecial subgroup as defined above. If $q^m = 3^1$, $E=M_{27}$, the group of $3\times3$ unitriangular matrices with elements in $\mathbb{F}_3$.

\end{Lemma}

\begin{proof}
Assume for contradiction that $q^m=3^1$ and $E=\mathbb{Z}_3 \rtimes \mathbb{Z}_9$. Then the elements of order $3$ form a characteristic subgroup $E_3$ of $E$. In fact, $E_3= \mathbb{Z}_3 \times \mathbb{Z}_3$ is abelian and non-cyclic. Now, $E_3$ is normal in $G_0$ as $E_3$ is characteristic in $E$ and $E$ is normal in $G_0$. This contradicts Theorem \ref{thm21} as $E_3$ is a normal, abelian subgroup of $G_0$ which is non-cyclic.
\end{proof}

Now we start at the bottom of the subgroup chain $Z \leq U \leq F \leq A \leq G_0 \leq G$ from Theorem \ref{thm21}. $Z=\mathbf{Z}(E) \cong \mathbb{Z}_q$ is a cyclic group of order $q$. $E$ is the extraspecial group of order $q^{2m+1}$ described above. $U$ is a cyclic group of order dividing $p^k-1$ (such that $p\neq q$) which acts irreducibly on a subspace $W$ of $V$ of dimension $k$. Further, $F=EU$ as a central product and $U=\mathbf{Z}(F)$. Then, $A/F$ acts faithfully on $E/Z$ so $A/F \lesssim SL(2m, p)$. Moreover, $G_0/A$ is isomorphic to a subgroup of Gal$(\mathbb{F}_{p^k} : \mathbb{F}_{p})$ which is a cyclic group of order $k$. So, $G_0/A \cong \mathbb{Z}_k$ or $G_0 = A$ as $k$ is prime. Finally, $G$ is the semidirect product of $G_0$ and a vector space of size $n:=p^d$.

\end{section}

\begin{section}{Computations}
In this section, we describe how we constructed candidates for groups $G_0$ with parameters, and checked each candidate to determine whether it had the desired properties, such as low rank, solvability, and primitivity. We carried out these computations in M\textsc{agma} and G\textsc{ap}.

\subsection{Irreducible Case} $G_0 \in \mathfrak{B}_I$. We use the following procedure to compute $G_0$ in the irreducible case.

{\bf Step 1:} Construct $E,$ the extraspecial subgroup of $G_0$ guaranteed by Theorem \ref{thm21}, as a subgroup of the general linear group $\mathrm{GL}(q^m,p^k)$ and its normalizer $N_E$ in $\mathrm{GL}(q^m,p^k).$

{\bf Step 2:} If $k>1$, embed $N_E$ into $\mathrm{GL}(kq^m,p)$ using the natural embedding $\mathrm{GL}(q^m,p^k) \to \mathrm{GL}(kq^m,p).$ If $k=1$, this step is unnecessary as $\mathrm{GL}(q^m,p^k) = \mathrm{GL}(kq^m,p).$

{\bf Step 3:} Construct $N := \mathbf{N}_{\mathrm{GL}(kq^m,p)}(N_E).$ Recall that $G_0$ will be a subgroup of $N$. Again, if $k=1$, this step is unnecessary as $\mathbf{N}_{\mathrm{GL}(kq^m,p)}(N_E)=N_E.$

{\bf Step 4:} Compute all subgroups of $N$ up to conjugacy class, and check each to determine if it is primitive and solvable, and moreover that it has rank less than five. Print the parameters of the groups that satisfy these conditions.

\subsection{Reducible Case} $G_0 \in \mathfrak{B}_R$. We modify the procedure slightly in the reducible case.

Since $d>kq^m$, after constructing $N_E$, we must embed it from $\mathrm{GL}(q^m,p^k)$ into $\mathrm{GL}(d/k,p^k)$ by taking its tensor product with $\mathrm{GL}(d/(kq^m),p^k)$. Then we proceed as in the irreducible case, replacing $kq^m$ with $d$ where appropriate. In all of the reducible cases we tested, $q^m=2^1, k=1$ so the embedding was simply from $\mathrm{GL}(2,p)$ to $\mathrm{GL}(d,p)$ with $d=4,6$ or $10$, and we didn't need to carry out the second embedding described in step 2 of the irreducible case.

\begin{Remark}
Step 4 is slow when the order of $N$ is large (i.e. mainly in the reducible cases). In these cases, we used M\textsc{agma}'s MaximalSubgroups command iteratively, discarding groups which were imprimitive or did not have low rank as their subgroups cannot be primitive or have low rank, respectively.
\end{Remark}

\end{section}

\begin{section}{Table of Results}

Below is a list of all groups in the $\mathfrak{B}_I$ class whose rank is less than or equal to four. For each case, we give the order of the maximal group with these properties, as well as its rank. We also provide the number of groups in each case that are subgroups of the maximal subgroup listed in the table. In the cases where $q=2$, we have handled the $E^+$ and $E^-$ cases separately and made a note of that. If only one of $E^+$ or $E^-$ is present for a given case, this is because the excluded case produced no valid groups.

We provide, in separate G\textsc{ap} \cite{GAP} files, all finite, primitive, solvable, linear groups of low rank. These groups are constructed as matrix groups acting on suitable fields, and each file corresponds to a single row in Table \ref{tbl41} or Table \ref{tbl42}.\footnote{\url{https://github.com/Yuxuan-Sun/Classifying-Solvable-Primitive-Permutation-Groups-of-Low-Rank}}

\begin{table}[H] 
\caption{Irreducible Cases}
\centering
\begin{tabular}{|L|L|L|L|L|L|L|L|L|L|}
\hline
\text{No.} & q & m & p & k & d & \operatorname{Rank}(G) & \operatorname{Max}|G_0| & \text{Num Gps}  & \text{Note} \\ \hline \hline
  1 & 2 & 1 & 3 & 1 & 2 & 2 & 48 & 4 & E- \\
  2 & 2 & 1 & 5 & 1 & 2 & 2 & 96 & 3 & E- \\
  3 & 2 & 1 & 7 & 1 & 2 & 2 & 144 & 7 & E- \\
  4 & 2 & 1 & 7 & 1 & 2 & 2 & 96 & 3 & E+ \\
  5 & 2 & 1 & 11 & 1 & 2 & 2 & 240 & 4 & E-\\
  6 & 2 & 1 & 23 & 1 & 2 & 2 & 528 & 3 & E-\\
  7 & 2 & 2 & 3 & 1 & 4 & 2 & 640 & 3 & E- \\\hline
  8 & 3 & 1 & 2 & 2 & 3 & 3 & 1296 & 7 & \\
  9 & 2 & 1 & 3 & 2 & 4 & 3 & 384 & 11 & E-\\
  10 & 2 & 1 & 13 & 1 & 2 & 3 & 288 & 2 & E-\\
  11 & 2 & 1 & 17 & 1 & 2 & 3 & 384 & 3 & E-\\
  12 & 2 & 1 & 19 & 1 & 2 & 3 & 432 & 3 & E-\\
  13 & 2 & 1 & 23 & 1 & 2 & 3 & 352 & 1 & E+ \\
  14 & 2 & 1 & 3 & 3 & 6 & 3 & 1872 & 6 & E-\\
  15 & 2 & 1 & 29 & 1 & 2 & 3 & 672 & 2 & E-\\
  16 & 2 & 1 & 31 & 1 & 2 & 3 & 720 & 1 & E-\\
  17 & 2 & 1 & 47 & 1 & 2 & 3 & 1104 & 1 & E-\\
  18 & 2 & 2 & 3 & 1 & 4 & 3 & 2304 & 13 & E+ \\
  19 & 2 & 2 & 7 & 1 & 4 & 3 & 1920 & 1 & E- \\ \hline
  20 & 3 & 1 & 7 & 1 & 3 & 4 & 1296 & 3 & \\
  21 & 2 & 1 & 5 & 2 & 4 & 4 & 1152 & 4 & E-\\
  22 & 2 & 1 & 31 & 1 & 2 & 4 & 480 & 1 & E+ \\
  23 & 2 & 1 & 37 & 1 & 2 & 4 & 864 & 1 & E-\\
  24 & 2 & 1 & 41 & 1 & 2 & 4 & 960 & 1 & E-\\
  25 & 2 & 1 & 43 & 1 & 2 & 4 & 1008 & 1 & E-\\
  26 & 2 & 1 & 53 & 1 & 2 & 4 & 1248 & 1 & E-\\
  27 & 2 & 1 & 59 & 1 & 2 & 4 & 1392 & 1 & E-\\
  28 & 2 & 1 & 71 & 1 & 2 & 4 & 1680 & 1 & E-\\
  29 & 2 & 2 & 5 & 1 & 4 & 4 & 4608 & 5 & E+ \\ \hline
\end{tabular}
\label{tbl41}
\end{table}

The second table lists all groups in the $\mathfrak{B}_R$ class whose rank $\le 4$. The formatting of this table follows that of the irreducible table. Note that we only found one case where the reducible case group $G_0$ was actually a maximal primitive linear group with rank $\leq4$. This group is displayed in line $30$ of the table below. In the other cases we checked, the group was either isomorphic to a subgroup of an irreducible case or had rank $>4$. There are two reasons why a group can be in both the  $\mathfrak{B}_R$ and $\mathfrak{B}_I$ classes.

Firstly, when $q=2$ and $k>1$ in Table \ref{tbl41}, the extraspecial group may act irreducibly over $\mathrm{GL}(q^m, p^k)$ but reducibly over $\mathrm{GL}(kq^m, p)$. Thus, despite the resulting $G_0$ groups being isomorphic, they belong to different classes. We note that in this case, the natural embedding $\mathrm{GL}(q^m,p^k)\to \mathrm{GL}(kq^m,p)$ provides a bijection between the two groups. In these cases, despite the vector spaces being represented differently, the difference only surfaces when we examine the vector space properties of $V$ since they are over different fields so the multiplication properties of the field are different. When treated purely as abelian groups, they are the same so the semidirect product is also isomorphic. The group action is the same and the semidirect product $G=V \rtimes G_0$ is isomorphic and we get one permutation group $G$ with a corresponding $G_0$ which can belong to either $\mathfrak{B}_I$ or $\mathfrak{B}_R$.

Secondly, due to the method which we used to construct these groups, starting with the extraspecial group and taking the normalizer in a general linear group, there are cases where two different extraspecial groups can lead to normalizers which contain the same group $G_0$. In this case, one of the extraspecial groups acts irreducibly and the other reducibly and both are contained in the final group $G_0$ which can thus belong to both $\mathfrak{B}_R$ and $\mathfrak{B}_I$.

\begin{table}[H]
\caption{Reducible Cases}
\centering
\begin{tabular}{|L|L|L|L|L|L|L|L|L|L|}
\hline
\text{No.} & q & m & p & k & d & \operatorname{Rank}(G) & \operatorname{Max}|G_0| & \text{Num Gps} & \text{Note} \\ \hline \hline
30 & 2 & 1 & 3 & 1 & 10 & 4 & 29040 & 1 & E-\\ \hline
\end{tabular}
\label{tbl42}
\end{table}
\begin{Remark}
Let $G_0$ be the group corresponding to line $30$ in the table above. For this final case where $|V|=p^{kd}=3^{10}$, we note that there is no corresponding irreducible case with $p^k=3^5, d=10$. This is because $G_0\leq\mathrm{GL}(10,3)$ is not entirely in the image of the natural embedding from $\mathrm{GL}(2,3^5)\to \mathrm{GL}(10,3)$. Thus, there is no subgroup of $\mathrm{GL}(2,3^5)$ isomorphic to $G_0$. What would have been the corresponding irreducible case is thus isomorphic to a proper subgroup of $G_0$ and we can then check that it has rank greater than $4$.
\end{Remark}

\end{section}

\begin{section}{Acknowledgements}

This research was conducted under NSF-REU grant DMS-1757233, DMS-2150205 and NSA grant H98230-21-1-0333, H98230-22-1-0022 by Dolorfino, Martin, Slonim, and Sun during the Summer of 2022 under the supervision of Yang. The authors gratefully acknowledge the financial support of NSF and NSA, and also thank Texas State University for providing a great working environment and support. Yang was also partially supported by grants from the Simons Foundation (\#499532, \#918096, YY).

We would also like to thank Professor Derek Holt for his invaluable help in understanding the computational side of this manuscript.

\end{section}

\end{document}